\documentclass[a4paper]{amsart}
\date{October 29, 2008}
\usepackage[dvips]{graphicx}
\usepackage{verbatim,enumerate}
\usepackage{amssymb}
\usepackage{times}

\usepackage{amsthm}
\theoremstyle{plain}
 \newtheorem{theorem}{Theorem}[section]
 \newtheorem*{theorem*}{Theorem}
 \newtheorem*{corollary*}{Corollary}
 \newtheorem*{lemma*}{Lemma}
 \newtheorem{proposition}[theorem]{Proposition}
 \newtheorem{lemma}[theorem]{Lemma}
 
 \newtheorem{introtheorem}{Theorem}
 \newtheorem{introcorollary}[introtheorem]{Corollary}
\theoremstyle{remark}
 
 \newtheorem{remark}[theorem]{Remark}
 \newtheorem*{remark*}{Remark}

\numberwithin{equation}{section}

\newcommand{\op}[1]{{\operatorname{ #1}}}

\newcommand{\SL}{\op{SL}}
\renewcommand{\O}{\op{O}}

\newcommand{\Per}{\op{Per}}
\newcommand{\rank}{\op{rank}}

\newcommand{\G}{\mathcal{G}}
\newcommand{\F}{\mathcal{F}}
\renewcommand{\S}{\mathcal{S}}
\renewcommand{\H}{\mathcal{H}}

\renewcommand{\Re}{\operatorname{Re}}
\renewcommand{\Im}{\operatorname{Im}}

\newcommand{\D}{\mathbb{D}}
\newcommand{\R}{\mathbb{R}}
\newcommand{\C}{\mathbb{C}}

\renewcommand{\L}{\mathbb{L}}

\newcommand{\imag}{\mathrm{i}}
\renewcommand{\epsilon}{\varepsilon}

\newcommand{\vect}[1]{\boldsymbol{#1}}

\newcommand{\trans}[1]{{\vphantom{#1}}^t#1}
\newcommand{\Hhol}{H^1_{\op{hol}}}

\title[Complete bounded curves]{
  Complete bounded holomorphic curves immersed in $\C^2$\\
  with arbitrary genus
}
\author{Francisco Martin}
\address[Martin]{
  Departamento de Geometr\'\i{}a y Topolog\'\i{}a,
  Universidad de Granada,
  18071 Granada, Spain.
}
\email{fmartin@ugr.es}
\author{Masaaki Umehara}
\address[Umehara]{%
  Department of Mathematics, Graduate School of Science,
  Osaka University,
  Toyonaka, Osaka 560-0043,
  Japan
}
\email{umehara@math.sci.osaka-u.ac.jp}
\author{Kotaro Yamada}
\address[Yamada]{%
  Faculty of Mathematics,
  Kyushu University,
  Fukuoka 812-8581, Japan%
}
\email{kotaro@math.kyushu-u.ac.jp}
\begin{document}
\begin{abstract}
In \cite{MUY}, a complete holomorphic 
immersion of the unit 
disk $\D$ into $\C^2$ whose image is bounded 
was constructed.
In this paper, we shall prove existence of 
complete holomorphic null
immersions of Riemann surfaces with arbitrary genus and finite
topology, whose image is bounded in $\C^2$.
To construct such immersions, we apply the method in \cite{Lopez}
to perturb the genus zero example in \cite{MUY} changing its genus. 

As an analogue the above construction,
we also give a new method to construct 
complete bounded minimal immersions
(resp.\ weakly complete maximal surface) 
with arbitrary genus and finite topology in Euclidean 3-space
(resp.\ Lorentz-Minkowski 3-space\-time).
\end{abstract}

\maketitle

\section{Introduction}

In \cite{MUY}, the authors constructed a complete bounded
minimal immersion of the unit disk $\D$ into $\R^3$ whose conjugate
is also bounded.
As applications of our results, we show in this paper that the technique
developed by F. J. L\'opez in \cite{Lopez} can be suitably  modified 
to give the following three things:
\begin{itemize}
\item The first examples of complete bounded complex submanifolds 
      with arbitrary genus immersed in $\C^2$;
\item A new and simple
      method to construct complete, bounded minimal surfaces 
      with arbitrary genus in the Euclidean 3-space;
\item A method to construct weakly complete, bounded maximal surfaces 
      with arbitrary genus in the Lorentz-Minkowski 3-spacetime.
\end{itemize}
Actually, a complete and bounded conformal minimal immersion of $\D$ into
$\R^3$ whose conjugate is also bounded is the real part of a complete and
bounded null holomorphic immersion
\[
   X_0:\D\longrightarrow \C^3.
\]

As it was shown in \cite{MUY}, if we consider the projection 
\begin{equation}\label{eq:proj2}
   \pi:\C^3\ni (z_1,z_2,z_3)\longmapsto (z_1,z_2)\in \C^2,
\end{equation}
the map $\pi\circ X_0:\D\to \C^2$ gives a complete bounded
complex submanifold immersed in $\C^2$.
By a perturbation of $X_0$, considering a similar deformation like
in \cite{Lopez}, we shall prove the following
\begin{introtheorem}\label{thm:main}
 Let $M$ be a closed Riemann surface.
 For an arbitrary positive integer $N$, 
 there exist points $Q_0$, $Q_1$, \dots, $Q_e$ $(e\ge N)$ on $M$, 
 neighborhoods $U_{l}$ of $Q_{l}$ which are homeomorphic to
 the unit disc $\D$, and a holomorphic immersion
 \[
    X:M\setminus \bigcup_{l=0}^e \overline{U_{l}}
               \longrightarrow \C^2
 \]
 which is complete and bounded.
\end{introtheorem}
If we do not specify the choice of $M$ with the fixed arbitrary genus, 
we can choose $M$ so that the number of ends is at most $3$
(see \cite[Lemma 3]{Lopez}).

Among other things, L\'opez in \cite{Lopez} gave a method to 
construct complete minimal surfaces between two
parallel planes of arbitrary genus. He constructed these examples
  as a perturbation of the minimal disk given by Jorge and Xavier
in \cite{jorge-xavier}.
However, the global boundedness was not treated in \cite{Lopez}, because
L\'opez' technique is not of use when you apply it directly to
Nadirashvili's minimal disk \cite{nadi}. 
So, the aim of this paper is to introduce some simple but 
appropriate changes in L\'opez'
machinery in order to obtain complete bounded minimal surfaces with
nontrivial topology.
Also the deformation parameter given in \cite{Lopez} is not sufficient to kill the
periods of the holomorphic immersions in $\C^2$, and so we must
consider another deformation space and calculate the new periods.

On the other hand, the proof of the above theorem allow us to conclude
that
\begin{introcorollary}\label{cor:cor}
 Let $M$ be a closed Riemann surface.
 For an arbitrary positive integer $N$, 
 there exist points $Q_0$,\dots,$Q_e$ $(e\ge N)$  on $M$,
 neighborhoods $U_l$ of $Q_l$ which are homeomorphic to
 the unit disc $\D$, and a complete conformal
 minimal immersion  into the  Euclidean $3$-space $\R^3$
 {\rm (}resp.\
 a weakly complete maximal surface in the sense of \cite{uy3},
 which may have singular points,
 whose first fundamental form is conformal
 on the regular set,
 in the Lorentz-Minkowski $3$-spacetime $\L^3${\rm)}
 \[
     x:M\setminus \bigcup_{l=0}^e \overline{U_l}\longrightarrow
     \R^3\quad(\text{resp.\ $\L^3$})
 \]
 whose image is bounded.
\end{introcorollary}

Though general methods to produce bounded complete minimal 
immersions into $\R^3$
with higher genus are already known
(cf.\ \cite{AFM}, \cite{FMM}, \cite{LMM-Handles} and \cite{LMM-Nonorientable} ),
this corollary gives a new and short way to construct such examples.  
On the other hand, for a maximal surface in $\L^3$,
this corollary produces the first examples
of weakly complete bounded maximal surfaces with arbitrary genus.
The {\em weak completeness\/} of maximal surfaces, which may admit
certain kind of singularities was defined in \cite{uy3}.
In \cite{alarcon} a weakly complete disk satisfying
a certain kind of boundedness was given. 
After that, a bounded example of genus zero with one end was
shown in \cite{MUY}.

Finally, we remark that the technique in this paper does not 
produce bounded null holomorphic curves immersed in $\C^3$
with arbitrary genus because of the difficulty in
obtaining the required
deformation parameters. 
If we do succeed to find an example, the above three examples could be all
realized as projections of it. 
So it is still an open problem to show the existence of complete bounded
null holomorphic immersion into $\C^3$ with arbitrary genus. 
\section{Preliminaries}\label{sec:prelim}
Let $M$ be a compact Riemann surface of genus $\kappa$, 
and fix a point $Q_0\in M$.
Let 
\begin{equation}\label{eq:gap}
   1=d_1 < d_2 < \dots < d_{\kappa} \leq 2\kappa-1
\end{equation}
be the {\em Weierstrass gap series\/} at $Q_0$, that is,
there exists a meromorphic function on $M$ which is holomorphic 
on $M\setminus\{Q_0\}$ such that $Q_0$ is the pole of order $d$
if and only if $d\not\in \{d_1,\dots,d_{\kappa}\}$.
Then there exists a meromorphic function $f$ on $M$ which is
holomorphic on $M\setminus\{Q_0\}$ and $Q_0$ is the pole of order
$m_0$ ($> d_\kappa$).
Let $\{Q_1,\dots,Q_e\}(\subset M)$ be the set of branch points of $f$.
Then the divisor $(df)$ of $df$ is written as 
\begin{equation}\label{eq:df}
  (df) = \frac{
            \prod_{l=1}^e Q_l^{m_l}
         }{Q_0^{m_0+1}}\qquad (m_0>d_{\kappa}),
\end{equation}
where $m_l$ ($l=1,\dots,e$) are positive integers, and the divisor
is expressed by the multiplication of these branch points.

\begin{remark}\label{rmk:N}
 In the latter construction of surfaces in Theorem~\ref{thm:main} 
 and Corollary~\ref{cor:cor},
 the ends of surfaces in $\C^2$ or $\R^3$ correspond to
 the points $Q_0$, \dots, $Q_e$. 
 We set
 \[
    f_{j+1}:=(f_j-c_j)^2,\quad f_0=f \qquad (j=1,2,3,\dots)
 \]
 inductively, where the $c_j$'s are complex numbers.
 Then, by replacing $f$ by
 $f_1$, $f_2$, $f_3$,\dots, we can make an example 
 such that $e$ is greater than the given number $N$.
\end{remark}
We write
\begin{equation}\label{eq:M0}
    M_0 := M \setminus \{Q_0,\dots,Q_e\}.
\end{equation}
Denoting by $\Hhol(M)$ and $\Hhol(M_0)$ the 
(first) holomorphic de Rham
cohomology group of $M$ and $M_0$
(as the vector space of holomorphic
differentials on $M$ factored by the subspace
of exact holomorphic differentials), respectively, then
\[
     \dim \Hhol(M) = \kappa,\qquad
     \dim \Hhol(M_0) = n 
\]
hold, where we set, for the sake of simplicity,
\begin{equation}\label{eq:n}
     n:=2\kappa + e.
\end{equation}
Take a basis $\{\zeta_1,\dots,\zeta_{\kappa}\}$ of $\Hhol (M)$.

\begin{lemma}[cf.\ III.5.13 in \cite{Kra}]
\label{lem:basis}
 One can take a basis of $\Hhol(M_0)$
 \[
    \{\zeta_1,\dots,\zeta_{\kappa};\xi_1,\dots,\xi_{\kappa};
      \eta_1,\dots,\eta_e\},
 \]
 where 
 $\xi_j$ is a meromorphic $1$-form on $M$ which is holomorphic on
 $M\setminus\{Q_0\}$,
 and $Q_0$ is a pole of order $m_j+1$,
 and
 $\eta_l$ is a meromorphic $1$-form on $M$ which 
 is holomorphic on
 $M\setminus\{Q_0,Q_l\}$,
 and $Q_0$ and $Q_l$  are poles of order $1$.
\end{lemma}

For simplicity, we set
\begin{equation}\label{eq:zeta-gen}
 \begin{aligned}
  \zeta_{\kappa+j} &:= \xi_j\qquad &(j=1,\dots,\kappa), \\
  \zeta_{2\kappa+l} &:= \eta_l  \qquad &(l=1,\dots,e).
 \end{aligned}
\end{equation}
Then $\{\zeta_1,\dots,\zeta_{n}\}$ ($n=2\kappa+e$)
is a basis of $\Hhol(M_0)$.

\begin{lemma}[{\cite[Lemma 1]{Lopez}}]\label{lem:v}
 There exists a meromorphic function $v$ on $M$ with the following
 properties{\rm :}
 \begin{enumerate}
  \item $v$ is holomorphic on $M_0=M\setminus\{Q_0,\dots,Q_e\}$,
  \item $Q_l$ $(l=1,\dots,e)$ is a pole of $v$ whose order 
	is greater than or equal to $m_l+2$, and
  \item $Q_0$ is a pole of $v$ whose order is greater than $m_0$.
 \end{enumerate}
\end{lemma}
\begin{proof}
 For each $l=1,\dots,e$, let $v_l$ be a meromorphic function on $M$
 which is holomorphic on $M\setminus\{Q_l\}$, and so that
 $Q_l$ is a pole of order $\max\{m_l+2,2\kappa\}$.
 On the other hand, take a meromorphic function $v_0$ on $M$
 which is holomorphic on $M\setminus\{Q_0\}$ and 
 $Q_l$ is a pole of order $\max\{m_0+1,2\kappa\}$.
 Then $v= v_0 + v_1 +\dots  + v_e$
 is a desired one.
\end{proof}

Using $f$ as in \eqref{eq:df} and $v$ as in Lemma~\ref{lem:v}, 
we define
\begin{equation}\label{eq:G-lambda}
  \G_{\Lambda}:= 
         \lambda_0 v + 
           \frac1{df}
             \sum_{j=1}^{n} \lambda_j \zeta_j
           : M \longrightarrow \C\cup\{\infty\},
\end{equation}
where
\[
    \Lambda:=(\lambda_0,\lambda_1,\dots,\lambda_{n})\in \C^{n+1}.
\]
\begin{lemma}[{\cite[Section 3]{Lopez}}]\label{lem:G-lambda}
 The function $\G_{\Lambda}$ as in \eqref{eq:G-lambda} is a 
 meromorphic function on $M$ such that
 \begin{enumerate}
  \item\label{item:G:1}
        $\G_{\Lambda}$ is holomorphic on
	$M_0=M\setminus\{Q_0,\dots,Q_e\}$,
  \item\label{item:G:2}
        if $\Lambda\neq \vect{0}$, $\G_{\Lambda}$ is 
	nonconstant on $M$, and
  \item\label{item:G:3}
        if $\lambda_0\neq 0$, $\G_{\Lambda}$ has poles
	at $Q_0,\dots,Q_e$.
 \end{enumerate}
\end{lemma}
\begin{proof}
 \ref{item:G:1} and \ref{item:G:2} are  obvious. 
 \ref{item:G:3} follows from the fact that $v$ has
 a  pole of higher order than $\zeta_1/df$, \dots, $\zeta_{n}/df$
 at $Q_j$ for each $j=0,\dots,e$.
\end{proof}
We write
\[
    |\Lambda|=\sqrt{|\lambda_0|^2+|\lambda_1|^2+\dots+|\lambda_{n}|^2},
\]
and consider the unit sphere in the space of $\Lambda$:
\[
  \S_1:=\{\Lambda\in \C^{n+1}\,;\, |\Lambda|=1\}.
\]
The following assertion is a modification of  
\cite[Lemma 2]{Lopez}, which is much easier to prove.
For our purpose, this weaker assertion is sufficient.

\begin{proposition}\label{prop:inverse}
 Let $\Lambda_0=(a_0,a_1,\dots,a_{n})$
 be a point in $\S_1$ satisfying $a_0\ne 0$.
 Then there exist  $\varepsilon(>0)$ and a neighborhood $W$ of
 $\Lambda_0$  in $\S_1$ such that if  $0<|t|<\varepsilon$ and
 $\Lambda\in W$, 
 the set $\G_{t\Lambda}^{-1}(\D)$  is conformally equivalent
 to a compact surface of genus $\kappa$ minus  
 $e+1$ pairwise disjoint discs with analytic regular boundaries.
 In particular, there are no branch points of $\G_{t\Lambda}$ on
 the boundary
 $\partial\G^{-1}_{t\Lambda}(\D)=\G^{-1}_{t\Lambda}(\partial \D)$.
\end{proposition}
\begin{proof}
 Since the poles of $\G_{\Lambda_0}$ are exactly
 $Q_0$, \dots, $Q_e$ and  
 $\G_{t\Lambda}^{-1}(\D)=\G_{\Lambda}^{-1}\bigl((1/t)\D\bigr)$,
 for sufficiently small $t$, the inverse image $\G_{t\Lambda}^{-1}(\D)$
 of the unit disk $\D$ by $\G_{t\Lambda}$ is homeomorphic to
 a compact surface of genus $\kappa$ minus  
 $e+1$ pairwise disjoint discs with piecewise analytic
 boundaries.
 Moreover, since the set of branch points of $\G_{\Lambda_0}$ does not
 have any accumulation points, $\G_{t\Lambda_0}^{-1}(\partial \D)$
 has no  branch points for sufficiently small $t$, and
 $\G_{t\Lambda}^{-1}(\partial \D)$ consists of real analytic regular
 curves in $M$.

 Furthermore, since $a_0\ne 0$, $\G_{\Lambda_0}/\G_{\Lambda}$
 are holomorphic near $Q_0$, \dots, $Q_e$ for any $\Lambda$
 which is sufficiently close to $\Lambda_0$.
 Thus $\G_{t\Lambda}^{-1}(\D)$
 has the same properties as $\G_{t\Lambda_0}^{-1}(\D)$.
\end{proof}

Under the situation in Proposition \ref{prop:inverse}, 
let $\{\gamma_1, \dots, \gamma_{2\kappa}\}$ be a family of loops on 
$\G_{\Lambda}^{-1}(\D)$ which is a homology basis of $M$.
On the other hand,
take a loop $\gamma_{2\kappa+l}$ on $\G_{\Lambda}^{-1}(\D)$
for each $l=1,\dots,e$, surrounding a neighborhood of $Q_l$
(as in Proposition~\ref{prop:inverse}).

We define the {\em period matrix\/} as
\begin{equation}\label{eq:period}
    P=(p_{kj}),\qquad
    p_{kj}:=\int_{\gamma_k}\zeta_j,
\end{equation}  
which is a nondegenerate $n\times n$ matrix.

\section{Proof of the main theorem} 
In this section, we give a proof of Theorem~\ref{thm:main} in the
introduction.

\subsection*{The initial immersion}
Let $X_0\colon{}\D\to \C^3$ be a complete holomorphic null immersion
whose image is bounded in $\C^3$ (as in Theorem~A in \cite{MUY}),
where $\D\subset \C$ is the unit disk.
We write
\begin{equation}\label{eq:initial}
    X_0(z) = \int_{0}^{z}
           \bigl( \varphi_1(z),\varphi_2(z),\varphi_3(z) \bigr)\,dz,
\end{equation}
where $z$ is a canonical coordinate on $\D\subset \C$ and
$\varphi_j$ ($j=1,2,3$) are holomorphic functions on $\D$.
Since $X_0$ is null, it holds that
\begin{equation}\label{eq:null}
   \bigl(\varphi_1\bigr)^2+
   \bigl(\varphi_2\bigr)^2+
   \bigl(\varphi_3\bigr)^2=0.
\end{equation}
Let $(g,\omega:=\omega_0\,dz)$ be the Weierstrass data of $X_0$, 
that is,
\begin{equation}\label{eq:w-data}
    \varphi_1 = (1-g^2)\omega_0,\quad
    \varphi_2 = \imag(1+g^2)\omega_0,\quad
    \varphi_3 = 2 g \omega_0,
\end{equation}
where $\imag=\sqrt{-1}$.
\newpage

\begin{lemma}\label{lem:phi}
 Let $X_0\colon{}\D\to \C^3$ be a null holomorphic immersion
 as in \eqref{eq:initial} whose image is not contained in any plane. 
 Then there exists a point $z_0\in \D$
 and a complex orthogonal transformation $T:\C^3\to \C^3$ 
 in $\O(3,\C)=\{a\in \op{GL}(3,\C)\,;\, \trans{a}=a^{-1}\}$
 {\rm (}$\trans{a}$ is the transpose of $a${\rm)}
 such that, up to replacing $X_0$ by $T\circ X_0$, 
 the following properties
 hold{\rm :}
 \begin{enumerate}
  \item\label{item:phi:1}
        $\varphi_1(z_0)=0$,
  \item\label{item:phi:2}
        $\varphi_3(z_0)\ne 0$ and $\varphi'_3(z_0)\ne 0$, where $'=d/dz$,
  \item\label{item:phi:3} 
        $\varphi_2(z_0)=\imag \varphi_3(z_0)$ and
	$\varphi'_2(z_0)=-\imag \varphi'_3(z_0)$.
 \end{enumerate}
 Moreover, if $X_0$ is complete and bounded, then so is $T\circ X_0$.
\end{lemma}
We shall now assume our initial $X_0$ satisfies the three properties
above, and set $z_0=0$ by a coordinate change of $\D$.

\begin{proof}[Proof of Lemma~\ref{lem:phi}]
 Since the image of $X_0$ is not contained in any plane,
 at least one of $\varphi_1$, $\varphi_2$ and $\varphi_3$,
 say $\varphi_3$, is not constant.
 Then we can take $z_0\in \D$ such that $\varphi_3(z_0)\ne 0$
 and $\varphi'_3(z_0)\ne 0$.
 Moreover, if $\varphi_1+\imag \varphi_2$ or $\varphi_1-\imag \varphi_2$ 
 vanishes identically, this contradicts that 
 the image of $X_0$ is not contained in any plane.
 So we may also assume that
 \[
     \varphi_1(z_0)\ne \pm \imag \varphi_2(z_0).
 \]
 When $\varphi_1(z_0)\ne 0$, 
 we replace $X_0$ by $T\circ X_0$, 
 where $T$ is the linear map
 associated with a complex orthogonal matrix
\[
    \begin{pmatrix}
       -c(1+c^2)^{-1/2} & \hphantom{-c}(1+c^2)^{-1/2} & 0 \\ 

       -       \hphantom{c}(1+c^2)^{-1/2} & -c(1+c^2)^{-1/2} & 0 \\ 
        0 & 0 & 1
    \end{pmatrix}
 \qquad (c=\varphi_2(z_0)/\varphi_1(z_0)).
\]
 Since $T$ is orthogonal, 
 $T\circ X_0$ is also a null holomorphic immersion.
 So we get the property \ref{item:phi:1}.
 Then \eqref{eq:null} implies that $\varphi_2(z_0)=\pm\imag\varphi_3(z_0)$.
 Replacing $\varphi_2$ by $-\varphi_2$ if necessary,  
 we may assume $\varphi_2(z_0)  = \imag \varphi_3(z_0)$.
 Differentiating \eqref{eq:null}, we have
 $\varphi_1(z_0)\varphi_1'(z_0)+ \varphi_2(z_0)\varphi_2'(z_0)+ 
  \varphi_3(z_0)\varphi_3'(z_0)=0$.
 Then 
 \[
   -\varphi_3(z_0)\varphi'_3(z_0)=\varphi_1(z_0)
             \varphi_1'(z_0) + 
       \varphi_2(z_0)\varphi_2'(z_0) 
           = \imag\varphi_3(z_0)\varphi_2'(z_0)
 \]
 holds.
 In particular, it holds that
 $\varphi_2'(z_0) = -\imag \varphi_3'(z_0)$.

 Next we prove the boundedness and completeness of $T\circ X_0$
 under the assumption that $X_0$ is complete and bounded.
 Since $X_0$ is bounded and $T$ is continuous, 
 $X_0(\D)$ and $T(X_0(\D))$ are both contained in the
 closed ball $\overline{B_0}(R)$ 
 in $\C^3$ of a certain radius $R>0$ centered at
 the origin. 
 We denote by $h_0$ the canonical metric on $\C^3$, and
 consider the pull-back metric $h_1=T^*h_0$ by $T$.
 Now we apply the following Lemma 3.1 in \cite{MUY}
 for $K:=\overline{B_0}(R)$ in $N:=\C^3$.
 Then there exist positive numbers $a$ and $b$ ($0<a<b$)
 such that $ah_0 < h_1 < bh_0$ on $\overline{B_0}(R)$.
 Now, we consider the pull back of the metric $h_0$ and
 $h_1$ by the immersion $X_0$. 
 Since
 \[
    (X_0)^*h_1=(X_0)^*(T^*h_0)
                   =(T\circ X_0)^*h_0=(X_1)^*h_0,
 \]
 we have that
 \[
         a (X_0)^*h_0 < (X_1)^*h_0 < b(X_0)^*h_0
 \]
 on $\D$. 
 Since $X_0$ is complete, $(X_0)^*h_0$ is a
 complete Riemannian metric on $\D$. 
 Then the above relation
 implies that $(X_1)^*h_0$ is also a complete Riemannian metric on
 $\D$,
 that is, $X_1$ is also a complete immersion.
\end{proof}

\subsection*{A family of holomorphic immersions}
Let 
\[
   \Lambda=(\lambda_0,\dots,\lambda_{n})\in\C^{n+1},\qquad
   \Delta=(\delta_1,\dots,\delta_{n})\in\C^{n},
\]
and take a meromorphic function $\G_{\Lambda}$ as in
\eqref{eq:G-lambda}.
Define a meromorphic function $\F_{\Delta}$ as
\begin{equation}\label{eq:F-delta}
    \F_{\Delta} =
    \frac{1}{df}\sum_{j=1}^{n}\delta_j\zeta_j
 \colon{}M\longrightarrow\C\cup\{\infty\}.
\end{equation}
and Weierstrass data $(\hat g,\hat\omega)$ on
$\G_{\Lambda}^{-1}(\D)$ as
\begin{equation}\label{eq:def-w-data}
 \hat g=\hat g_{(\Lambda,\Delta)}:= h_{\Delta}\cdot (g\circ \G_{\Lambda}),\quad
 \hat\omega =\hat\omega_{(\Lambda,\Delta)}:=
     \frac{\omega_0\circ \G_{\Lambda}\,df}{h_{\Delta}} \quad
     \bigl(h_{\Delta}:=\exp \F_{\Delta}\bigr),
\end{equation}
where $f$ is the meromorphic function as in \eqref{eq:df},
and define holomorphic $1$-forms on $\G_{\Lambda}^{-1}(\D)$ as
\begin{equation}\label{eq:mod-forms}
  \Psi_1 = (1-\hat g^2)\hat\omega,\qquad
  \Psi_2 = \imag(1+\hat g^2)\hat\omega,\qquad
  \Psi_3 = 2\hat g \hat\omega.
\end{equation}

The following lemma is a modified version of \cite[Theorem 3]{Lopez}.
(In fact, our data \eqref{eq:F-delta} and \eqref{eq:def-w-data} for the 
surfaces are somewhat different from those in \cite{Lopez}.)
\begin{lemma}\label{lem:complete}
 If $X_0$ as in \eqref{eq:initial} is a complete immersion, 
 the metric
 \[
      d\hat s^2:=(1+|\hat g|^2)^2 |\hat\omega|^2
 \]
 determined by 
 $(\hat g,\hat\omega)$ as in \eqref{eq:def-w-data} is a complete
 Riemann metric  on $\G_{\Lambda}^{-1}(\D)$  for a sufficiently small 
 $(\Lambda,\Delta)\neq (\vect{0},\vect{0})$.
\end{lemma}
\begin{proof}
 As in the equations (15) and (17) in \cite{Lopez}, there exists a positive
 constant
 $a$ ($<1$) such that 
 \[
   a<|h_{\Delta}|<\frac1a \qquad\text{and}\qquad
   \left|\frac{df}{d\G_{\Lambda}}\right|>a\qquad
   \text{on $\G_{\Lambda}^{-1}(\D)$}.
 \]
 Then, (setting $z=\G_{\Lambda}$)
 \begin{align*}
    (1+|\hat g|^2)|\hat\omega|
    &=
    (1+|gh_{\Delta}|^2)\left|\frac{\omega_0}{h_{\Delta}}\right|
       \left|\frac{df}{d\G_{\Lambda}}\right|\,|dz|\\
    &\geq (a^2+|ag|^2)|a\omega_0|\,(a\,|dz|)=a^4(1+|g|^2)|\omega|.
 \end{align*}
 Thus we have the conclusion.
\end{proof}
Thus, for each (sufficiently small) 
$(\Lambda,\Delta)\in \C^{2n+1}\setminus\{\vect{0}\}$,
there exists a complete null immersion
\begin{equation}\label{eq:F-deform}
   X_{(\Lambda,\Delta)}:=
    \int_{z_0}^z
      \bigl(\Psi_1,\Psi_2,\Psi_3\bigr)
      \colon{}\widetilde{\G_{\Lambda}^{-1}(\D)}\longrightarrow
   \C^3,\end{equation}
where $\widetilde{\G_{\Lambda}^{-1}(\D)}$ denotes the universal cover of
$\G_{\Lambda}^{-1}(\D)$. 
In fact, the line integral $\int_{z_0}^z \Psi_j$ ($j=1,2,3$) 
from a base point $z_0$ depends on the choice of the path,
but can be considered as a 
single-valued function on $\widetilde{\G_{\Lambda}^{-1}(\D)}$.

Then we get the following assertion, which can be proven exactly
in the same way as Corollary B in \cite{MUY}:
\begin{proposition}\label{cor:complete-proj}
 Let $\pi$ be the projection as in \eqref{eq:proj2}.
 Then $\pi\circ X_{(\Lambda,\Delta)}$ is 
 a complete immersion  of $\widetilde{\G_{\Lambda}^{-1}(\D)}$ 
 into $\C^2$.
\end{proposition}

\subsection*{The period map}
Under the situations above, we define the {\em period map}
\begin{multline}\label{eq:period-map}
 \Per_1\colon{}\C^{2n+1}\ni(\Lambda,\Delta)\\
     \longmapsto 
     \trans{\left(
      \int_{\gamma_1}\Psi_1,\dots,
      \int_{\gamma_{n}}\Psi_1,
      \int_{\gamma_{1}}\Psi_2,\dots,
      \int_{\gamma_{n}}\Psi_2
     \right)}\in
     \C^{2n},
\end{multline}
where $n=2\kappa+e$ (see \eqref{eq:n}),
``$\trans{(~)}$'' is the transposing operation for matrices,
$\gamma_j$'s are loops as in \eqref{eq:period}, and
$\Psi_1$ and $\Psi_2$ are as in \eqref{eq:mod-forms}.
The following assertion is an analogue of \cite[Theorem 2]{Lopez}:

\begin{proposition}\label{prop:jacobian}
Suppose that $X_0$ satisfies the three conditions as in
Lemma \ref{lem:phi}. Then
the $(2n)\times (2n)$ matrix
 \begin{equation}\label{eq:jacobian}
    J_1 := \left.
       \left(
	 \frac{\partial \Per_1}{\partial\lambda_1},\dots,
	 \frac{\partial \Per_1}{\partial\lambda_{n}},
	 \frac{\partial \Per_1}{\partial\delta_1},\dots,
	 \frac{\partial \Per_1}{\partial\delta_{n}}
       \right)
       \right|_
      {(\Lambda,\Delta)=(\vect{0},\vect{0})}
 \end{equation}
 is nondegenerate.
\end{proposition}
\begin{proof}
 Note that
 \begin{equation}\label{eq:FG-initial}
    \left. \G_{\Lambda}\right|_{\Lambda=\vect{0}}
        =\left. \F_{\Delta}\right|_{\Delta=\vect{0}}=0,\qquad
    \left. h_{\Delta}\right|_{\Delta=\vect{0}}=1.
 \end{equation}
 By the definitions, we have
 \[
    \left.\frac{\partial \G_{\Lambda}}{\partial \lambda_j}\right|_
    {\Lambda=\vect{0}}
    = \frac{\zeta_j}{df},\qquad
    \left.\frac{\partial h_{\Delta}}{\partial \delta_j}\right|_
    {\Delta=\vect{0}}=
    \left.\frac{\partial \exp \F_{\Delta}}{\partial \delta_j}\right|_
    {\Delta=\vect{0}}= \frac{\zeta_j}{df}
 \]
 for $j=1,\dots,n$.
 Then
 \begin{align*}
   \left.\frac{\partial \hat g}{\partial\lambda_j}
      \right|_{(\Lambda,\Delta)=(\vect{0},\vect{0})} &
    = \left(\left.\frac{dg}{dz}\right|_{z=0}\right)
      \left(\left.\frac{\partial \G_{\Lambda}}{\partial\lambda_j}
      \right|_{\Lambda=\vect{0}}\right)
    = g'(0)\frac{\zeta_j}{df},\\
   \left.\frac{\partial \hat \omega}{\partial\lambda_j}
      \right|_{(\Lambda,\Delta)=(\vect{0},\vect{0})}
   &  = \omega_0'(0)\zeta_j,\\
   \left.\frac{\partial \hat g}{\partial\delta_j}
      \right|_{(\Lambda,\Delta)=(\vect{0},\vect{0})}
   & = g(0)\frac{\zeta_j}{df},\qquad
   \left.\frac{\partial \hat \omega}{\partial\delta_j}
      \right|_{(\Lambda,\Delta)=(\vect{0},\vect{0})}
     =-\omega_0(0)\zeta_j,
 \end{align*}
 where $'=d/dz$.
 Hence we have
 \begin{align*}
  \left.\frac{\partial\Psi_1}{\partial\lambda_j}
  \right|_{(\Lambda,\Delta)=(\vect{0},\vect{0})}
  &=\left.
  \left(-2\hat g\frac{\partial \hat g}{\partial \lambda_j}\hat\omega+
        (1-\hat g^2)\frac{\partial\hat\omega}{\partial\lambda_j}\right)
  \right|_{(\Lambda,\Delta)=(\vect{0},\vect{0})}\\
  &= \bigl(-2 g(0)g'(0)\omega_0(0)+(1-g(0)^2)\omega_0'(0)\bigr)\zeta_j\\
  &= \left.\bigl((1-g^2)\omega\bigr)'\right|_{z=0}\zeta_j
  = \varphi_1'(0)\zeta_j,\\
  \left.\frac{\partial\Psi_1}{\partial\delta_j}
  \right|_{(\Lambda,\Delta)=(\vect{0},\vect{0})}
  &=  -(1+\{g(0)\}^2)\omega_0(0)=\imag\varphi_2(0)\zeta_j,
\end{align*}
where $\varphi_j$'s are holomorphic functions on $\D$ as in \eqref{eq:w-data}.
Similarly, we have
 \[
  \frac{\partial\Psi_2}{\partial\lambda_j}
  = \varphi'_2(0)\zeta_j,\quad
  \frac{\partial\Psi_3}{\partial\lambda_j}
  = \varphi'_3(0)\zeta_j,\quad
  \frac{\partial\Psi_2}{\partial\delta_j}
  = -\imag\varphi_1(0)\zeta_j,\quad
  \frac{\partial\Psi_3}{\partial\delta_j}
  =0
 \]
at $(\Lambda,\Delta)=(\vect{0},\vect{0})$.
 Thus, we have that
 \begin{equation}\label{eq:periods2}
 \begin{alignedat}{2}
   \frac{\partial}{\partial\lambda_j}\int_{\gamma_k}\Psi_1
     &=\varphi_1'(0)\int_{\gamma_k}\zeta_j,\qquad
  &\frac{\partial}{\partial\delta_j}\int_{\gamma_k}\Psi_1
     &=\imag\varphi_2(0)\int_{\gamma_k}\zeta_j,\\
  \frac{\partial}{\partial\lambda_j}\int_{\gamma_k}\Psi_2
     &=\varphi_2'(0)\int_{\gamma_k}\zeta_j,\qquad
  &
   \frac{\partial}{\partial\delta_j}\int_{\gamma_k}\Psi_2
     &=-\imag\varphi_1(0)\int_{\gamma_k}\zeta_j,\\
  \frac{\partial}{\partial\lambda_j}\int_{\gamma_k}\Psi_3
     &=\varphi_3'(0)\int_{\gamma_k}\zeta_j,\qquad
  &
   \frac{\partial}{\partial\delta_j}\int_{\gamma_k}\Psi_2
     &=0
 \end{alignedat}
 \end{equation}
 hold at $(\Lambda,\Delta)=(\vect{0},\vect{0})$,
 for $j,k=1,\dots,n$.
 Hence the matrix $J_1$ in \eqref{eq:jacobian} is written as
 \[
    J_1 = \begin{pmatrix}
	  \varphi_1'(0) P & \hphantom{-}\imag\varphi_2(0) P \\
          \varphi_2'(0) P & -\imag\varphi_1(0) P
	  \end{pmatrix}=
          \begin{pmatrix}
	  \varphi_1'(0) P & \hphantom{-}\imag\varphi_2(0) P \\
          \varphi_2'(0) P & O
	  \end{pmatrix}
          \qquad (\varphi_2(0), \varphi'_2(0)\neq 0)
 \]
 because of Lemma~\ref{lem:phi},
 where $P$ is the nondegenerate period matrix as in \eqref{eq:period}.
 Hence $J_1$ is nondegenerate.
\end{proof}

\subsection*{The period-killing problem}
Since $\Per_1(\vect{0})=\vect{0}$,
Proposition~\ref{prop:jacobian} yields that 
there exists a holomorphic map
$c \mapsto \bigl(\lambda_1(c),\dots,\lambda_{n}(c),
             \delta_1(c),\dots,\delta_{n}(c)
	     \bigr)$
such that
\[
   \Per_1\bigl(c,\lambda_1(c),\dots,\lambda_{n}(c),
         \delta_1(c),\dots,\delta_{n}(c)
      \bigr)  = 0 
\]
for sufficient small $c$.
We set
\[
   \G_c = \G_{\Lambda(c)},\qquad
   \text{where}\quad
   \Lambda(c):=\bigl(c,\lambda_1(c),\dots,\lambda_{n}(c)\bigr).
\]
Since 
$\Lambda(0)=\vect{0}$,
there exists an analytic function $\Lambda^*(c)$ such that
$\Lambda(c)=c\Lambda^*(c)$ near $c=0$.
Now we can apply Proposition \ref{prop:inverse} to
$\Lambda_0:=\Lambda^*(0)/|\Lambda^*(0)|\in\S_1$.
Then, for sufficiently small $c$, $\G_{c}^{-1}(\D)$ 
is conformally equivalent to $M$ minus $e+1$ 
pairwise disjoint discs with analytic regular boundaries, 
and the map
\begin{multline}\label{eq:result}
   X_c:=\pi\circ X_{\bigl(\Lambda(c),\Delta(c)\bigr)},\\
   \Lambda(c)=(c,\lambda_1(c),\dots,\lambda_{n}(c)),\quad
   \Delta(c)=(\delta_1(c),\dots,\delta_{n}(c))
\end{multline}
is well-defined on $\G_{c}^{-1}(\D)$.
Moreover, by Corollary~\ref{cor:complete-proj}, 
$X_c$ is a complete immersion for any sufficient small $c$.

\subsection*{Boundedness of $X_c$}
By Proposition~\ref{prop:inverse}, $d\G_c$ does not vanish
on $\partial\G_c^{-1}(\D)$ for sufficiently small $c$.
Then if we choose
a real number $r\in (0,1)$ sufficiently close to $1$,
we have
\[
    d\G_c\neq 0 \qquad \text{on}\quad \G_c^{-1}(\overline\D\setminus\D_r),
\]
where 
$\D_r=\{z\in\C\,;\,|z|<r\}$.
Moreover, $\G_c^{-1}(\overline\D\setminus\D_r)$ is exactly
a union of  $e+1$ closed annular domains surrounding
the points $Q_0,Q_1,\dots,Q_e$ in $M$.
To show the boundedness of $X_c$, it is sufficient to show that
the image of each annular domain by $X_c$ is  bounded.
For the sake of simplicity, we show the boundedness
of $X_c$ at $Q_0$. 
We denote by $\overline\Omega$ the closed
annular domain surrounding the point $Q_0$.
Then
\[
   \G:=\G_c|_{\overline\Omega}:\overline\Omega
    \longrightarrow \overline\D\setminus\D_r
\]
gives a holomorphic finite covering.

\begin{figure}
\begin{center}
 \includegraphics{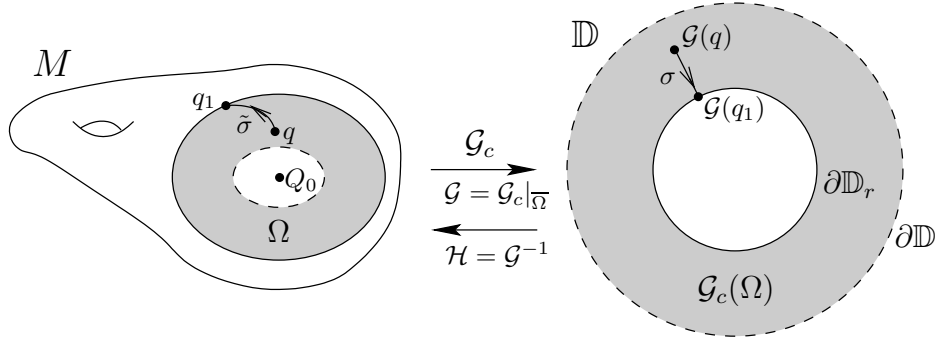}
\end{center}
\caption{The segments $\tilde\sigma$ and $\sigma$}%
\label{fig:sigma}
\end{figure}

Since $\G_c^{-1}(\overline\D_r)\subset \G_c^{-1}(\D)$ is compact and 
$X_c\colon{}\G_c^{-1}(\D)\to\C^2$ is holomorphic, there exists
a positive constant $K_0$ such that
\begin{equation}\label{eq:bound1}
   |X_c| \leq K_0\qquad \text{on}\quad \G_c^{-1}(\overline\D_r).
\end{equation}
We denote by $\Omega$ the set of interior points of $\overline\Omega$,
and fix $q\in \Omega$ arbitrarily.
Let $\sigma(t)$ ($0\le t\le 1$)
be a line segment such that
$\sigma(0)=\G(q)$ and $\sigma(1)=r \G(q)/|\G(q)|\in \partial \D_r$, 
that is,
\[
   \sigma(t):=(1-t)\G(q)+t\frac{r\G(q)}{|\G(q)|}.
\]
Since $\G$ is a covering map, there exists a unique smooth curve 
$\tilde \sigma:[0,1]\to \overline\Omega$
such that $\tilde \sigma(0)=q$ and $\G\circ \tilde\sigma=\sigma$.
Moreover, there exists a neighborhood $U$ of the line segment
$\sigma([0,1])$ and a holomorphic map
$\H:U\to M$ which gives the (local) inverse of $\G$.
By definition, we have
$\tilde \sigma=\H\circ \sigma$.
We set 
\[
   q_1:=\tilde \sigma(1)\in 
        \partial \Omega \setminus \G_c^{-1}(\partial \D),
\]
that is, $q_1$ lies on the connected component of $\partial\Omega$
which is further from $Q_0$, see Figure~\ref{fig:sigma}.

By \eqref{eq:bound1}, it is sufficient to show that
\begin{equation}\label{eq:small-part}
  \bigl|
    X_c(q_1)-X_c(q)
  \bigr|
  =\left|\int_{\tilde \sigma}(\Psi_1,\Psi_2)\right|
\end{equation}
is bounded from above by a constant which does not depend on $q$.
Here we have
\begin{multline}\label{eq:int-phi-one}
 \left|\int_{\tilde \sigma}\Psi_1\right|=
 \left|
   \int_{\tilde \sigma}
         \left(1-\bigl(h_{\Delta}\bigr)^2 
                 \bigl(g\circ \G_c\bigr)^2\right)
          \frac{\omega_0\circ \G_c}{%
                h_{\Delta}}df
 \right|\\
\begin{aligned}
 &=
 \left|
   \int_0^1
    \left(
        \left(1-\bigl(g\circ \sigma(t)\bigr)^2\bigl(h(t)\bigr)^2\right) 
        \frac{\omega_0\circ \sigma(t)}{h(t)}
        \left. \frac{d\bigl(f\circ \H(z)\bigr)}{dz}\right|_{z=\sigma(t)}
        \frac{d\sigma}{dt}
    \right)\,dt
 \right|\\
  &=|I_1+I_2|,
\end{aligned}
\end{multline}
where $h(t):=h_{\Delta}\circ \tilde \sigma(t)$ and
\begin{align*}
   I_1 &= \frac{1}{2}\int_0^1
           \left(
            \left. \bigl(\varphi_1(z)-\imag\varphi_2(z)\bigr) 
           \right|_{z=\sigma(t)}
            \frac{d\sigma}{dt}
           \right)
           \left(
             \left. \frac{1}{h_{\Delta}\circ \H(z)}
             \frac{df\circ \H(z)}{dz}\right|_{z=\sigma(t)}  
           \right)\,dt,\\
   I_2 &= \frac{-1}{2}\int_0^1
           \left(
            \left. \bigl(\varphi_1(z)+\imag\varphi_2(z)\bigr) 
           \right|_{z=\sigma(t)}
            \frac{d\sigma}{dt}
           \right)
           \left(
             \left. h_{\Delta}\circ \H(z)
             \frac{df\circ \H(z)}{dz}\right|_{z=\sigma(t)}  
           \right)\,dt.
\end{align*}
Here we used the relations $2\omega_0=\varphi_1-\imag\varphi_2$ and 
$2g^2\omega_0=-(\varphi_1+\imag\varphi_2)$.
To estimate $I_1$, 
we shall apply Lemma~\ref{lem:app} in the appendix for
\[
   a(t) =
            \left. \bigl(\varphi_1(z)-\imag\varphi_2(z)\bigr) 
\right|_{z=\sigma(t)}
            \frac{d\sigma}{dt}
	    \quad\text{and}\quad
   b(t)	=  
             \left. \frac{1}{h_{\Delta}\circ \H(z)}
             \frac{df\circ \H(z)}{dz}\right|_{z=\sigma(t)}.  
           \]
Let us check the hypotheses:
\[
   (A(s):=)\int_0^s a(t)\,dt = \int_{\sigma([0,s])} 
      (\varphi_1(z)-\imag\varphi_2(z))\,dz
      = \left.\bigl(X_1(z)-\imag X_2(z)\bigr)\right|_{z=\sigma(0)}^{z=\sigma(s)},
\]
where $X=X_0$ as in \eqref{eq:initial}.
Since $X$ is bounded, $|A(s)|$ is bounded for all $s$.
On the other hand, 
\[
    b(t) = \tilde b\bigl(\sigma(t)\bigr),\qquad\text{where}\quad
    \tilde b(z) =
       \frac{1}{h_{\Delta}\circ \H(z)}\frac{d(f\circ \H(z))}{dz}.
\]
Since $\H(z)$ can be considered as a single valued holomorphic
function on a certain finite covering of $\overline\D\setminus \D_r$,
both $\tilde b(z)$ and $\tilde b'(z)$ ($'=d/dz$) are 
bounded by a constant. 
Hence by the lemma, we have that $|I_1|$ is bounded.
Similarly, we can show that $|I_2|$ is bounded, and
we can conclude that
the integration of $\Psi_1$ 
(and similarly $\Psi_2$) along $\tilde\sigma$ is bounded.

The resulting immersion $X_c$ has an arbitrary number of
ends by setting $f=f_N$ as in Remark \ref{rmk:N}.

\section{Proof of the corollary}
In this section we shall prove Corollary~\ref{cor:cor} in the
introduction.
Define $\G_{\Lambda}$, $(\Psi_1,\Psi_2,\Psi_3)$ and 
the null immersion
$X_{(\Lambda,\Delta)}$ as in the previous section,
and define real parameters $(s_j,t_j)$ as
\begin{equation}\label{eq:real-parameter}
    \lambda_j = s_j + \imag t_j,\qquad
    \delta_j  = s_{n+j}+ \imag t_{n+j}.
\end{equation}
Note that for a holomorphic function $F(u)$ in $u=s+\imag t$, one has:
\begin{equation}\label{eq:re-im}
\begin{alignedat}{2}
    \frac{\partial\Re F}{\partial s} &=\Re \frac{dF}{du},\qquad
   &\frac{\partial\Re F}{\partial t} &= -\Im \frac{dF}{du},\\
    \frac{\partial\Im F}{\partial s} &=\Im \frac{dF}{du},\qquad
   &\frac{\partial\Im F}{\partial t} &=\hphantom{-} \Re \frac{dF}{du}.
\end{alignedat}
\end{equation}
\subsection*{Minimal surfaces in $\R^3$}
First, we treat the case of minimal surfaces in $\R^3$.
Let
\[
   x= x_{(\Lambda,\Delta)}:=\Re X_{(\Lambda,\Delta)}=
       \Re\int(\Psi_1,\Psi_2,\Psi_3)\colon{}
       \widetilde{\G_{\Lambda}^{-1}(\D)}\longrightarrow \R^3,
\]
where $\widetilde{\G_{\Lambda}^{-1}(\D)}$ is the universal cover of
$\G_{\Lambda}^{-1}(\D)$.
Then $x$ is a conformal minimal immersion, and the induced metric
is complete because of Lemma~\ref{lem:complete}.
If $x_{(\Lambda,\Delta)}$ were well-defined on $\G_{\Lambda}^{-1}(\D)$
for sufficiently small $(\Lambda,\Delta)$, boundedness of
the image of $x_{(\Lambda,\Delta)}$ can be proved in a similar way 
as in the previous section.
Hence it is sufficient to solve the period-killing problem
to show the corollary.

We define the period map
\begin{multline}\label{eq:period-R}
 \Per_2\colon{}\R^{4n+2}
      \ni(s_0,\dots,s_{2n};
          t_0,\dots,t_{2n})
  \\
     \longmapsto 
     \trans{\left(
      \left(\Re\!\!\int_{\gamma_k}\Psi_1\right)_{k=1,\dots,n},
      \left(\Re\!\!\int_{\gamma_k}\Psi_2\right)_{k=1,\dots,n},
      \left(\Re\!\!\int_{\gamma_k}\Psi_3\right)_{k=1,\dots,n}
     \right)}\in
     \R^{3n},
\end{multline}
where $(s_j,t_j)$ are real parameter as in \eqref{eq:real-parameter}.
Consider the $(3n)\times (4n)$ matrix
\begin{equation}\label{eq:Jacobi-R}
 J_2:=\left.
  \left(
    \frac{\partial\Per_2}{\partial s_1},\dots,
    \frac{\partial\Per_2}{\partial s_{2n}},
    \frac{\partial\Per_2}{\partial t_1},\dots,
    \frac{\partial\Per_2}{\partial t_{2n}}
  \right)\right|_{(\Lambda,\Delta)=(\vect{0},\vect{0})}.
\end{equation}
To solve the period-killing problem, 
it is sufficient to show that the rank of $J_2$ is $3n$.

By \eqref{eq:periods2} and \eqref{eq:re-im}, we have
\[
  J_2 = \begin{pmatrix}
	 \Re \bigl(\varphi_1'(0) P\bigr) &
	 \hphantom{-}\Re \bigl(\imag\varphi_2(0) P\bigr) &
	 -\Im \bigl(\varphi_1'(0) P\bigr) &
	 -\Im \bigl(\imag\varphi_2(0) P \bigr)\\
	 \Re  \bigl(\varphi_2'(0) P\bigr) &
	 -\Re \bigl(\imag\varphi_1(0) P\bigr) &
	 -\Im \bigl(\varphi_2'(0) P\bigr) &
	 \hphantom{-}\Im  \bigl(\imag\varphi_1(0) P\bigr) \\
	 \Re  \bigl(\varphi_3'(0) P \bigr)&
	 O &
	 -\Im \bigl(\varphi_3'(0) P\bigr) &
         O
	\end{pmatrix},
\]
where the $n\times n$ matrix $P$ is the period matrix as in \eqref{eq:period}.
Here, we remark that the real vectors are linearly independent over $\R$
if and only if they are linearly independent over $\C$.
Since we may assume that $X_0$ satisfies the three conditions as in
Lemma~\ref{lem:phi}, we have that
\begin{align*}
   \rank J_2 &=
   \rank\begin{pmatrix}
	  \varphi_1'(0) P & 
	  \hphantom{-}\imag\varphi_2(0) P & 
	  \overline{\varphi_1'(0) P} & 
	  \hphantom{-}\overline{\imag\varphi_2(0)P}\\
	  \varphi_2'(0) P & 
	  -\imag\varphi_1(0) P & 	  
	  \overline{\varphi_2'(0) P} & 
	  -\overline{\imag\varphi_1(0) P} \\
	  \varphi_3'(0) P & O &
	  \overline{\varphi_3'(0) P} & O
	 \end{pmatrix}\\
       &=              
       \rank\begin{pmatrix}
	  \hphantom{-\imag}\varphi_1'(0) P & 
	  \imag\varphi_2(0) P & 
	  \hphantom{\imag}\overline{\varphi_1'(0) P} & 
	  \overline{\imag\varphi_2(0)P}\\
	  -\imag\varphi_3'(0) P & 
          O &
	  \imag\overline{\varphi_3'(0) P} & 
          O \\
	  \hphantom{-\imag}\varphi_3'(0) P & O &
	  \hphantom{\imag}\overline{\varphi_3'(0) P} & O
	 \end{pmatrix}\\
   &=\rank
     \begin{pmatrix}
       O & P & * & * \\
       P & O & * & * \\
       O & O & 2\overline P & *
     \end{pmatrix}= 3n,
\end{align*}
 since $P$ is nondegenerate. Then we can solve the period problem as we
 did in the previous section.

\subsection*{Maximal surfaces in the Lorentz-Minkowski spacetime}
We denote by $\L^3$ the Lorentz-Minkowski $3$-spacetime, that is,
$\bigl(\R^3;(x_0,x_1,x_2)\bigr)$ endowed with the indefinite metric
$-(dx_0)^2+(dx_1)^2+(dx_2)^2$.
Under the same settings as above, we set
\[
   y = y_{(\Lambda,\Delta)}
     := \Re\int (\imag \Psi_1,\Psi_2,\Psi_3)\colon{}
       \widetilde{\G_{\Lambda}^{-1}(\D)}\longrightarrow \L^3.
\]
Then $y$ gives a maximal surface (a mean curvature zero surface), 
which could possibly have singular points. 
In particular, since the holomorphic lift
\begin{equation}\label{eq:maxface-lift}
    \int (\imag \Psi_1,\Psi_2,\Psi_3)\colon{}
       \widetilde{\G_{\Lambda}^{-1}(\D)}\longrightarrow \C^3
\end{equation}
is an immersion, $y$ is a {\em maxface\/} in the sense of \cite{uy3}.
Moreover, the induced metric by \eqref{eq:maxface-lift} is complete
because of Lemma~\ref{lem:complete}.
Hence $y$ is a {\em weakly complete maxface\/} in the sense of
\cite{uy3}.

To show the $\L^3$ case of Corollary~\ref{cor:cor}, 
we consider the period map
\begin{multline}\label{eq:period-L}
 \Per_3\colon{}\R^{4n+2}
      \ni(s_0,\dots,s_{2n};
          t_0,\dots,t_{2n})
  \\
     \longmapsto 
     \trans{\left(
      \left(\Im\!\!\int_{\gamma_k}\Psi_1\right)_{k=1,\dots,n},
      \left(\Re\!\!\int_{\gamma_j}\Psi_2\right)_{k=1,\dots,n},
      \left(\Re\!\!\int_{\gamma_j}\Psi_3\right)_{k=1,\dots,n}
     \right)}\in
     \R^{3n}.
\end{multline}
Then 
\begin{equation}\label{eq:Jacobi-L}
 J_3:=\left.
  \left(
    \frac{\partial\Per_3}{\partial s_1},\dots,
    \frac{\partial\Per_3}{\partial s_{2n}},
    \frac{\partial\Per_3}{\partial t_1},\dots,
    \frac{\partial\Per_3}{\partial t_{2n}}
  \right)\right|_{(\Lambda,\Delta)=(\vect{0},\vect{0})}
\end{equation}
has the expression
\begin{align*}
  J_3 =\begin{pmatrix}
	 \Im \bigl(\varphi_1'(0) P\bigr) &
	 \hphantom{-}\Im \bigl(\imag\varphi_2(0) P\bigr) &
	 \hphantom{-}\Re \bigl(\varphi_1'(0) P\bigr) &
	 \Re \bigl(\imag\varphi_2(0) P \bigr)\\
	 \Re  \bigl(\varphi_2'(0) P\bigr) &
	 -\Re \bigl(\imag\varphi_1(0) P\bigr) &
	 -\Im \bigl(\varphi_2'(0) P\bigr) &
	 \Im  \bigl(\imag\varphi_1(0) P\bigr) \\
	 \Re  \bigl(\varphi_3'(0) P \bigr)&
	 O &
	 -\Im \bigl(\varphi_3'(0) P\bigr) &
         O
	\end{pmatrix},
\end{align*}
and then it can be easily checked that
 $J_3$ is of rank $3n$ like as in the case of $\R^3$.
Hence we conclude as in the previous cases.

\appendix
\section{A lemma to show boundedness}
In this appendix, we show the following lemma:
\begin{lemma}\label{lem:app}
 Let $a$ and $b\colon{}I\to \C$ be smooth functions, where 
 $I=[0,L]$ is an interval in $\R$.
 Suppose that there exist constants $C_1$, $C_2$ and $C_3$ such that
 \[
    \left|\int_0^s a(t)\,dt\right|< C_1,\qquad
    |b(s)|<C_2,\quad\text{and}\quad |b'(s)|<C_3
 \]
 for all $s\in I$.
 where $'=d/ds$.
 Then there exists a constant $C=C(C_1,C_2,C_3,L)$ such that
 \[
     \left|\int_I a(t)b(t)\,dt\right|< C.
 \]
\end{lemma}
\begin{proof}
 Let
 \[
    A(s)= \int_0^s a(t)\,dt.
 \]
 Then
 \begin{align*}
   \left|\int_I a(t)b(t)\,dt\right|
   &= 
   \left|\int_I A'(t)b(t)\,dt\right|
   =
   \left|\left. A(t)b(t)\right|^{t=L}_{t=0}
         -\int_I A(t)b'(t)\,dt
   \right|\\
   &=
      \left|A(L)b(L)  -\int_I A(t)b'(t)\,dt
   \right|
   \leq |A(L)b(L)|+\left|\int_I A(t)b'(t)\right|\\
   &\leq C_1C_2 + L C_1C_3. &\qed
 \end{align*}
\renewcommand{\qed}{\relax}
\end{proof}


\end{document}